\documentclass{article}%
\usepackage{amsfonts}
\usepackage{amsmath}
\usepackage{amssymb}
\usepackage{graphicx}%
\newtheorem{theorem}{Theorem}

\newtheorem{corollary}[theorem]{Corollary}

\newtheorem{definition}[theorem]{Definition}

\newtheorem{notation}[theorem]{Notation}

\newtheorem{proposition}[theorem]{Proposition}
\newtheorem{remark}[theorem]{Remark}

\newenvironment{proof}[1][Proof]{\noindent\textbf{#1.} }{\ \rule{0.5em}{0.5em}}
\begin{document}

\title{{Composition and exponential of compactly supported generalized
integral kernel operators }}
\author{S.Bernard \thanks{Universit\'{e} des Antilles et de la Guyane, Laboratoire
AOC, Campus de Fouillole, 97159 Pointe-\`{a}-Pitre, Guadeloupe, E-mail:
severine.bernard@univ-ag.fr}, J.-F.Colombeau \thanks{Institut Fourier,
Universit\'{e} J. Fourier, 100 rue des Maths, BP 74, 38402 St Martin d'Heres,
France, E-mail: jf.colombeau@wanadoo.fr}, A.Delcroix \thanks{Same address as
the first author, E-mail: antoine.delcroix@univ-ag.fr}}
\date{}
\maketitle

\begin{abstract}
We extend the theory of distributional kernel operators to a
framework of generalized functions, in which they are replaced by
integral kernel operators. Moreover, in contrast to the
distributional case, we show that these generalized integral
operators can be composed unrestrictedly. This leads to the
definition of the exponential of a subclass of such operators.

\end{abstract}


\section{Introduction}

The theory of nonlinear generalized functions
\cite{colombeau84,colombeau85,colombeau92,grosser,nedeljkov},
which appears as a natural extension of the theory of
distributions, seems to be a suitable framework to overcome the
limitations of the classical theory of unbounded operators.

Following a first approach done by D. Scarpalezos in
\cite{scarpalezos}, we introduced in \cite{BCD1} a natural concept
of integral kernel operators in this setting. In addition, we
showed that these operators are characterized by their kernel. Our
approach has some relationship with the one of \cite{garetto}, but
is less restrictive and uses other technics of proofs. Let us
quote that classical operators with smooth or distributional
kernel can be canonically extended in the framework of generalized
functions, through the sheaf embeddings of $\mathcal{C}^{\infty
}(\cdot)$ or $\mathcal{D}^{\prime}(\cdot)$ into
$\mathcal{G}(\cdot)$, the sheaf of spaces of generalized
functions. This shows that our theory is a natural extension of
the classical one.

After recalling briefly the mathematical framework, we focus on
the case of generalized integral operators with compactly
supported kernel.\ We show that such operators can be composed
unrestrictedly and that their composition is still a generalized
integral operator with a kernel having a compact support.

This allows to consider their iterate composition and the question
of summation of series of such operators naturally arise. It has
been solved for the exponential, with additional assumptions on
the growth of the kernel with respect to the scaling parameter, in
view of applications to theoretical physics.

Let us mention that the question of composition of generalized
integral kernel operators has been investigated in more general
cases in \cite{BCD1}, namely for integral operators with properly
supported kernel and with kernel in the algebra
$\mathcal{G}_{L^{2}}$, constructed from $H^{\infty}$.\ In this
last case, the exponential of generalized integral kernel
operators can also be defined with the above mentioned assumptions
on the kernel.

\section{Generalized integral operators}

\subsection{The sheaf of nonlinear generalized functions}

In this section, we recall briefly some elements of the theory of generalized
numbers and functions. We refer the reader to
\cite{biagioni,colombeau84,colombeau85,colombeau92,grosser,marti03,nedeljkov}
for more details.

Let $\mathrm{C}^{\infty}(\cdot)$ be the sheaf of complex valued
smooth functions on $\mathbb{R}^{d}$ ($d\in\mathbb{N}$) endowed
with the usual topology of uniform convergence of all the
derivatives on compact sets.\ For every open set $\Omega$ of
$\mathbb{R}^{d}$, this topology can be described by the family of
semi norms $\left( p_{K,l}(\cdot)\right)
_{K\Subset\Omega,l\in\mathbb{N}}$ with $p_{K,l}(f)=\sup_{x\in
K,\,\left\vert \alpha\right\vert \leq l}\left\vert
\partial^{\alpha}f\left(  x\right)  \right\vert $, for all $f$ in $\mathrm{C}^{\infty}(\Omega)$ (the notation $K\Subset
\Omega$ means that\ the set $K$ is a compact set included in
$\Omega $).\smallskip

Set, with\ $\mathcal{E}\left(  \Omega\right)  =\mathrm{C}^{\infty}\left(
\Omega\right)  ^{\left(  0,1\right]  }$,%
\begin{align*}
\mathcal{E}_{M}\left(  \Omega\right)   &  =\left\{  \left(  f_{\varepsilon
}\right)  \in\mathcal{E}\left(  \Omega\right)  \left\vert \forall
K\Subset\Omega,\ \forall l\in\mathbb{N},\;\exists q\in\mathbb{N}%
,\;p_{K,l}\left(  f_{\varepsilon}\right)  =\mathrm{O}\left(  \varepsilon
^{-q}\right)  \;\mathrm{as}\;\varepsilon\rightarrow0\right.  \right\}, \\
\mathcal{I}\left(  \Omega\right)   &  =\left\{  \left(  f_{\varepsilon
}\right)  \in\mathcal{E}\left(  \Omega\right)  \left\vert \forall
K\Subset\Omega,\ \forall l\in\mathbb{N},\;\forall p\in\mathbb{N}%
,\;p_{K,l}\left(  f_{\varepsilon}\right)  =\mathrm{O}\left(  \varepsilon
^{p}\right)  \;\mathrm{as}\;\varepsilon\rightarrow0\right.  \right\}  .
\end{align*}

The functor $\mathcal{E}_{M}:\Omega\rightarrow\mathcal{E}_{M}\left(
\Omega\right)  $ (\textit{resp.} $\mathcal{I}:\Omega\rightarrow\mathcal{I}%
\left(  \Omega\right)  $) defines a sheaf of subalgebras of the sheaf
$\mathcal{E}\left(  \cdot\right)  $ (\textit{resp}. a sheaf of ideals of the
sheaf $\mathcal{E}_{M}\left(  \cdot\right)  $) \cite{marti03}.

\begin{definition}
The sheaf of factor algebras $\mathcal{G}\left(  \cdot\right)  =\mathcal{E}%
_{M}\left(  \cdot\right)  /\mathcal{I}\left(  \cdot\right)  $ is called the
sheaf of \emph{Colombeau type algebras}.
\end{definition}

The sheaf $\mathcal{G}$ turns to be a sheaf of differential algebras and a
sheaf of modules on the factor ring $\overline{\mathbb{C}}=\mathcal{X}\left(
\mathbb{C}\right)  /\mathcal{N}\left(  \mathbb{C}\right)  $ with
\begin{align*}
\mathcal{X}\left(  \mathbb{K}\right)   &  =\left\{  \left(  r_{\varepsilon
}\right)  \in\mathbb{K}^{\left(  0,1\right]  }\,\left\vert \,\exists
q\in\mathbb{N},\;\left\vert r_{\varepsilon}\right\vert =\mathrm{O}\left(
\varepsilon^{-q}\right)  \;\mathrm{as}\;\varepsilon\rightarrow0\right.
\right\}, \\
\mathcal{N}\left(  \mathbb{K}\right)   &  =\left\{  \left(  r_{\varepsilon
}\right)  \in\mathbb{K}^{\left(  0,1\right]  }\,\left\vert \,\forall
p\in\mathbb{N},\;\left\vert r_{\varepsilon}\right\vert =\mathrm{O}\left(
\varepsilon^{p}\right)  \;\mathrm{as}\;\varepsilon\rightarrow0\right.
\right\}  ,
\end{align*}
where $\mathbb{K}=\mathbb{C}$ or $\mathbb{K}=\mathbb{R}$.

\begin{notation}
For $\left(  f_{\varepsilon}\right)  \in\mathcal{E}_{M}\left(  \Omega\right)
$, $Cl\left(  f_{\varepsilon}\right)  $ will denote its class in $\mathcal{G}%
\left(  \Omega\right)  $.
\end{notation}

As $\mathcal{G}$ is a sheaf, the notion of support of a section $f\in
\mathcal{G}\left(  \Omega\right)  $ ($\Omega$ open subset of $\mathbb{R}^{d}$)
makes sense. Thus the \emph{support} of a generalized function $f\in\mathcal{G}%
\left(  \Omega\right)  $, denoted by $supp\,f$, is the complement
in $\Omega$ of the largest open subset of $\Omega$ where $f$ is
null. We denote by $\mathcal{G}_{C}\left(  \Omega\right)  $ the
subset of elements of $\mathcal{G}\left( \Omega\right)$ with
compact supports. In particular, such compactly supported
generalized functions have the following property: Every
$f\in\mathcal{G}_{C}$ has a representative $\left(
f_{\varepsilon}\right) \in\mathcal{E}_{M}\left(  \Omega\right)  $
such that each $f_{\varepsilon}$ has the same compact support. We
say that such a representative has a \textit{global compact
support}.

\subsection{Definitions and first properties}

Let $X$ (\textit{resp}. $Y$) be an open subset of $\mathbb{R}^{m}$
(\textit{resp}. $\mathbb{R}^{n}$). We denote by $\mathcal{G}$$_{ps}(X\times
Y)$ the set of generalized functions $g$ of $\mathcal{G}$$(X\times Y)$
properly supported in the following sense:%
\begin{equation}%
\begin{array}
[c]{c}%
\forall~O_{1}\subset X\text{ relatively compact open subset}%
,\ \ \ \ \ \ \ \ \ \ \ \ \ \ \ \ \ \ \ \ \ \\
\ \ \ \ \ \ \ \ \ \ \ \ \ \ \ \ \ \ \ \ \ \ \ \ \ \ \ \ \ \exists K_{2}\Subset
Y~/~supp\,g\cap(O_{1}\times Y)~\subset O_{1}\times K_{2}.
\end{array}
\label{HypGPS}%
\end{equation}
The set $\mathcal{G}$$_{ps}(X\times Y)$ is clearly a subalgebra of
$\mathcal{G}$$(X\times Y)$.

\begin{proposition}
\label{intgps}\cite{BCD1}$~$For $g$ in $\mathcal{G}$$_{ps}(X\times
Y)$, there exists $G\in$ $\mathcal{G}$$(X)$ such that, for all
relatively compact open subset $O_{1}$ of $X$,
\[
G_{\left\vert O_{1}\right.  }=Cl\left(  \left(  x\mapsto\int_{K_{2}%
}g_{\varepsilon}(x,y)\,\mathrm{d}y\right)  _{\left\vert O_{1}\right.
}\right)
\]
where $\left(  g_{\varepsilon}\right)  $ is a representative of $g$ and
$K_{2}\Subset Y$ is such that $supp\,g\cap(O_{1}\times Y)~\subset O_{1}\times
K_{2}.$
\end{proposition}

\begin{notation}
With a slight abuse, we shall denote $G=\int
g(\cdot,y)\,\mathrm{d}y$ or $G\left(  \cdot_{1}\right)  =\int
g(\cdot_{1},y)\,\mathrm{d}y$.
\end{notation}

\begin{definition}
Let $H$ be in $\mathcal{G}$$_{ps}(X\times Y)$. We call \emph{generalized
integral operator} the map
\[
\widehat{H}~:~%
\begin{array}
[t]{cll}%
\mathcal{G}(Y) & \rightarrow & \mathcal{G}(X)\\
f & \mapsto & \widehat{H}(f)=\int H(\cdot,y)f(y)\,\mathrm{d}y.
\end{array}
\]
We say that $H$ is the \emph{kernel} of the generalized integral operator
$\widehat{H}$.
\end{definition}

This map is well defined, due to proposition \ref{intgps} since
the application $g=H(\cdot_{1},\cdot_{2})f\left( \cdot_{2}\right)
$ is in $\mathcal{G}$$_{ps}(X\times Y)$, for all
$f\in\mathcal{G}(Y)$. Moreover, it is linear.

\begin{remark}
\label{GF04RemKomp}Any $H\in$ $\mathcal{G}$$(X\times Y)$ compactly
supported satisfies (\ref{HypGPS}) and $\widehat{H}$ is well
defined.\ Moreover, a straightforward calculation shows that the
image of $\widehat{H}$ is included in $\mathcal{G}$$_{C}(X)$.
Furthermore, the definition of $\widehat{H}$ does not need to
refer to proposition \ref{intgps} in this case. Indeed, if $H$ is
in $\mathcal{G}_{C}(X\times Y)$ with $supp\,H\subset K_{1}\times
K_{2}$ ($K_{1}\Subset X$, $K_{2}\Subset Y$) and $f$ in
$\mathcal{G}(Y)$, we have
\[
\widehat{H}(f)=Cl\left(  x\mapsto\int_{K_{2}}H_{\varepsilon}%
(x,y)f_{\varepsilon}(y)\,\mathrm{d}y\right)
\]
where $\left(  H_{\varepsilon}\right)  $ (resp. $\left(  f_{\varepsilon
}\right)  $) is any representative of $H$ (resp. $f$).
\end{remark}

\begin{remark}
If $H$ is in $\mathcal{G}$$(X\times Y)$ without hypothesis on the support, we
can define a map $\widehat{H}~$: $\mathcal{G}$$_{C}(Y)\rightarrow$
$\mathcal{G}$$(X)$ in the same way, since for all $f$ in $\mathcal{G}$%
$_{C}(Y)$, the function $H(\cdot_{1},\cdot_{2})f(\cdot_{2})$ is in
$\mathcal{G}$$_{ps}(X\times Y)$. In this case, the generalized integral
operator could also be defined globally since $f$ has a representative with a
global compact support, as quoted above.
\end{remark}

This case leads us to make the link between the classical theory
of integral operators acting on $\mathcal{D}(Y)$ and the
generalized one. Indeed, if $h$ belongs to
$\mathcal{D}^{\prime}(X\times Y)$ and $\widehat{h}$ is the
classical operator of kernel $h$, then the following diagram is commutative:%
\[%
\begin{array}
[c]{ccc}%
\mathcal{D}(Y) & \overset{\widehat{h}}{\rightarrow} & \mathcal{D}^{\prime
}(X)\\
\downarrow\sigma &  & \downarrow i_{S}\\
\mathcal{G}_{C}(Y) &
\overset{\widehat{i_{S}^{\prime}(h)}}{\rightarrow} &
\mathcal{G}(X),
\end{array}
\]
where $\sigma$ (\textit{resp}. $i_{S}$, $i_{S}^{\prime}$) is the
usual embedding of $\mathcal{D}(Y)$ into $\mathcal{G}_{C}(Y)$
(\textit{resp}. $\mathcal{D}^{\prime}(X)$ into $\mathcal{G}(X)$~,
$\mathcal{D}^{\prime }(X\times Y)$ into $\mathcal{G}(X\times Y)$).
This shows that our theory extends \textquotedblleft
canonically\textquotedblright\ the classical one. We refer to
\cite{BCD1} for more details on the relationship with classical
cases and to \cite{delcroix2,nedeljkov} for the definition of the
sheaves embeddings of $\mathcal{D}\left(  \cdot\right)  $ and
$\mathcal{D}^{\prime }(\cdot)$ into $\mathcal{G}(\cdot)$.

\begin{remark}
\label{SKT}The map$~\widehat{}~:~$$\mathcal{G}$$_{ps}(X\times Y)\rightarrow
$$\mathcal{L}$$($$\mathcal{G}$$(Y),$$\mathcal{G}$$(X))$ is a linear map of
$\overline{\mathbb{C}}$-modules. Moreover, $\widehat{H}$ is
continuous for the sharp topologies \cite{scarpalezos}.\newline
Conversely, the third author showed in \cite{delcroix1} that any
continuous linear map from $\mathcal{G}_{C}(Y)$ to
$\mathcal{G}(X)$, satisfying appropriate growth hypothesis with
respect to the regularizing parameter $\varepsilon$, can be
written as a generalized integral kernel operator.
\end{remark}

The following result shows that the map $\widehat{}$ , defined in
remark \ref{SKT}, is injective.

\begin{theorem}
\label{caracker}\cite{BCD1, garetto2}~Characterization of
generalized integral operators by their kernel: One has
$\widehat{H}=0$ if and only if $H=0$.
\end{theorem}

\subsection{Composition of generalized integral operators}

For this topic, we only consider in this paper generalized
integral operators with compactly supported kernel.

\begin{theorem}
\label{corol}For $H$ in $\mathcal{G}_{C}(X\times\Xi)$ and $K$ in
$\mathcal{G}_{C}(\Xi\times Y)$, $\widehat{H}\circ\widehat{K}:\mathcal{G}%
(Y)\rightarrow\mathcal{G}_{C}(X)$ is a generalized integral operator whose
kernel $L$ is an element of $\mathcal{G}$$_{C}(X\times Y)$ defined globally by
$L(\cdot_{1},\cdot_{2})=\int_{\Xi}H(\cdot_{1},\xi)K(\xi,\cdot_{2}%
)\,\mathrm{d}\xi$.\newline Moreover, there exists $K_{1}$ (resp.
$K_{2}~$, $K_{3}$) a compact set of $X$ (resp. $\Xi~,Y$) such that
the support of $H$ (resp. $K$) is contained in the interior of
$K_{1}\times K_{2}$ (resp. $K_{2}\times K_{3}$). In this case, the
support of $L$ is contained in $K_{1}\times K_{3}$.
\end{theorem}

\begin{proof}
For all $f$ in $\mathcal{G}$$(Y)$, $\widehat{K}\left(  f\right)  $
is well defined and belongs to $\mathcal{G}_{C}$$(\Xi)$, according
to remark \ref{GF04RemKomp}.\ This allows the definition of the
composition $\widehat {H}\circ\widehat{K}$. Let us verify now the
assertion concerning the support of $L$. Since $H$ (\textit{resp}.
$K$) is in $\mathcal{G}_{C}(X\times\Xi)$ (\textit{resp}.
$\mathcal{G}_{C}(\Xi\times Y)$), we can find $K_{1}$
(\textit{resp.} $K_{2}~$, $K_{3}$) satisfying the second assertion
of the theorem. Then,
\[
L(\cdot_{1},\cdot_{2})=\int_{K_{2}}H(\cdot_{1},\xi)K(\xi,\cdot_{2}%
)\,\mathrm{d}\xi
\]
is a well defined generalized function, according to the theory of
integration of generalized functions on compact sets
\cite{colombeau84,colombeau92}. Denote by $\left(
H_{\varepsilon}\right)  $ (\textit{resp}. $\left(
K_{\varepsilon}\right)  $) a representative of $H$ (\textit{resp}.
$K$) and set $O_{1}=X\setminus K_{1}$, $O_{3}=Y\setminus K_{3}$.
The map
$L_{\varepsilon}(\cdot_{1},\cdot_{2})=\int_{K_{2}}H_{\varepsilon}(\cdot
_{1},\xi)K_{\varepsilon}(\xi,\cdot_{2})\,\mathrm{d}\xi$ is a
representative of $L$. For $U\Subset X$ and $V\Subset Y$ such that
$U\times V\subset O_{1}\times O_{3}$, we have either $U\subset
O_{1}$ or $V\subset O_{3}$. We shall suppose, for example, that
$U\subset O_{1}$. For $(x,y)\in U\times V$, we have \
\[
\left.  \left\vert L_{\varepsilon}(x,y)\right\vert =\left\vert \int_{K_{2}%
}H_{\varepsilon}(x,\xi)K_{\varepsilon}(\xi,y)\,\mathrm{d}\xi\right\vert
\right.  \leq Vol(K_{2})p_{U\times
K_{2}},_{0}(H_{\varepsilon})p_{K_{2}\times
V},_{0}(K_{\varepsilon}),
\]
where $Vol(K_{2})$ denotes the volume of $K_2$. Therefore
\begin{equation}
p_{U\times V},_{0}(L_{\varepsilon})\leq Vol(K_{2})p_{U\times K_{2}}%
,_{0}(H_{\varepsilon})p_{K_{2}\times V},_{0}(K_{\varepsilon}).
\label{GF04BCD11}%
\end{equation}
As $(H_{\varepsilon\left\vert O_{1}\times\Xi\right.  })$ is in $\mathcal{I}%
(O_{1}\times\Xi)$ and $U\cap K_{2}\subset O_{1}\times\Xi$, it
follows that $p_{U\times K_{2},0}(H_{\varepsilon})=O\left(
\varepsilon^{m}\right)  $ as $\varepsilon\rightarrow0$, for all
$m\in\mathbb{N}$. Moreover, $\left( K_{\varepsilon}\right)$ is in
$\mathcal{E}_{M}\left(  \Xi\times Y\right)  $. Thus, relation
(\ref{GF04BCD11}) implies that $p_{U\times V},_{0}(L_{\varepsilon
})=O\left(  \varepsilon^{m}\right)  $ as
$\varepsilon\rightarrow0$, for all $m\in\mathbb{N}$. Finally,
$(L_{\varepsilon})$ satisfies the null estimate of order $0$ for
all compact sets included in $O_{1}\times O_{3}$. Using theorem
1.2.3 of \cite{grosser}, we can conclude, without estimates on the
derivatives, that $L_{\left\vert O_{1}\times O_{3}\right. }=0$.
Therefore, the support of $L$ is contained in $K_{1}\times K_{3}$.
From this, a straightforward verification, using once more the
integration on compact sets of generalized functions, shows that
$\widehat{H}\circ\widehat{K}=\widehat{L}$. \medskip
\end{proof}

Repeted applications of theorem \ref{corol} show the following:

\begin{corollary}
\label{corocompo}For $H$ in $\mathcal{G}_{C}(X^{2})$ and all $n\geq2$, the
composition $n$ times of $\widehat{H}$ is a well defined operator $\widehat
{H}^{n}$ with image in $\mathcal{G}_{C}(X)$. Moreover, $\widehat{H}^{n}$admits
as kernel $L_{n}\in\mathcal{G}_{C}(X^{2})$ defined by%
\[
L_{n}(\cdot_{1},\cdot_{2})=\int_{X^{n-1}}H(\cdot_{1},\xi_{1})H(\xi_{1},\xi
_{2})\cdots H(\xi_{n-1},\cdot_{2})\,\mathrm{d}\xi_{1}\mathrm{d}\xi_{2}%
\cdots\mathrm{d}\xi_{n-1}.
\]
Furthermore, for all $n\geq2$, the support of $L_{n}$ is contained in the one
of $H$.
\end{corollary}

\section{Application: exponential of generalized integral operators}

In this section, we define the exponential of generalized integral operators
in a particular case and study some of their properties. We need before to
introduce a convenient subsheaf of $\mathcal{G}(\cdot)$. For $\Omega$ open set
of$\mathbb{\ R}^{d}$ ($d\in\mathbb{N}$), set
\[
\mathcal{H}_{ln}(\Omega)=\left\{  (u_{\varepsilon})%
\in\mathcal{E(}\Omega)~/~\forall K\Subset\Omega,~\forall l\in\mathbb{N}%
,~~p_{K,l}(u_{\varepsilon})=O(|\ln\varepsilon|)\text{ as }\varepsilon
\rightarrow0\right\}  .
\]
The set $\mathcal{H}_{ln}(\Omega)$ is a linear subspace of $\mathcal{E}%
_{M}(\Omega)$ (but not a subalgebra). Define%
\[
\mathcal{G}_{ln}(\Omega)=\frac{\mathcal{H}_{ln}(\Omega)}{\mathcal{I}(\Omega
)}\mbox{ and }
\mathcal{G}_{C\,ln}(\Omega)=\mathcal{G}_{ln}(\Omega)\cap
\mathcal{G}_{C}(\Omega).
\]

\begin{theorem}
\label{defexp1} Let $H$ be in $\mathcal{G}_{C\,ln}(X^{2})$. Denote
by $L_{n}$ the kernel of
$\widehat{H}^{n}:\mathcal{G}(X)\rightarrow \mathcal{G}_{C}(X)$
defined as in corollary \ref{corocompo} and $\left(
L_{n,\varepsilon}\right)  $ a representative of $L_{n}$. For all
$\varepsilon\in(0,1]$, the series
$\sum_{n\geq1}\frac{L_{n,\varepsilon}}{n!}$ (by setting $L_1=H$)
normally converges, for the usual topology of uniform convergence
on compact subsets of $X^{2}$.\ Denote by $S_{\varepsilon}$ its
sum.\ The net $\left( S_{\varepsilon}\right)  $ belongs to
$\mathcal{E}_{M}(X^{2})$. Furthermore, $S=Cl\left(
S_{\varepsilon}\right)  $ defines a compactly supported element of
$\mathcal{G}(X^{2})$ only depending on $H$.\newline The well
defined operator $e^{\widehat{H}}=\widehat{S}+Id$ (where $Id$ is
the operator identity) will be called the \emph{exponential} of
$\widehat{H}$.
\end{theorem}

The \textbf{proof} of this theorem can be divided in three parts.
The first part contains the estimates of
$\sum_{n\geq1}\frac{L_{n,\varepsilon}}{n!}$ for a particular
representative of $L_{n}$, given by a fixed representative of $H$.
The second part deals with the independence of $Cl\left(
S_{\varepsilon }\right)  $ with respect to the chosen
representative of $L_{n}$, that is of $H$. The third part shows
that $S$ is compactly supported.

We shall give here mainly the first part of this proof and refer the reader to
\cite{BCD1} for other parts. Let $H$ be in $\mathcal{G}_{C\,ln}(X^{2})$ and
$\left(  H_{\varepsilon}\right)  $ one of its representative. According to
corollary \ref{corocompo}, we have $\widehat{H}^{n}=\widehat{L_{n}%
}:\mathcal{G}(X)\rightarrow\mathcal{G}_{C}(X)$ and $L_{n}\in$ $\mathcal{G}%
_{C}(X^{2})$ admits as representative $\left(  L_{n,\varepsilon}\right)  $
with
\[
L_{n,\varepsilon}(\cdot_{1},\cdot_{2})=\int_{K^{n-1}}H_{\varepsilon}(\cdot
_{1},\xi_{1})H_{\varepsilon}(\xi_{1},\xi_{2})\cdots H_{\varepsilon}(\xi
_{n-1},\cdot_{2})\,\mathrm{d}\xi_{1}\mathrm{d}\xi_{2}\cdots\mathrm{d}\xi
_{n-1},
\]
where $K$ is a compact set of $X$ such that the support of $H$ is contained in
the interior of $K^{2}$.

For all compact subset of $X^{2}$ of the form $K_{1}\times K_{2}$, $\left(
\alpha,\beta\right)  \in\mathbb{N}^{d}\times\mathbb{N}^{d}$ and $(x,y)\in
K_{1}\times K_{2}$, one has%
\begin{align*}
\left\vert \frac{\partial^{\alpha+\beta}L_{2,\varepsilon}}{\partial x^{\alpha
}\partial y^{\beta}}(x,y)\right\vert  &  =\left\vert \int_{K}\frac
{\partial^{\alpha}H_{\varepsilon}}{\partial x^{\alpha}}(x,\xi)\frac
{\partial^{\beta}H_{\varepsilon}}{\partial y^{\beta}}(\xi,y)\,\mathrm{d}%
\xi\right\vert \\
&  \leq\int_{K}p_{K_{1}\times K,|\alpha|}(H_{\varepsilon})p_{K\times
K_{2},|\beta|}(H_{\varepsilon})\,\mathrm{d}\xi.
\end{align*}
It follows that
\[
p_{K_{1}\times K_{2},|(\alpha,\beta)|}(L_{2,\varepsilon})\leq Vol(K)p_{V^{2}%
,|(\alpha,\beta)|}^{2}(H_{\varepsilon}),
\]
where $V$ is a compact set of $X$ containing $K$,$K_{1}$ and $K_{2}$. By an
iterative method, we show that, for all $n\geq2$,
\[
p_{K_{1}\times K_{2},|(\alpha,\beta)|}(L_{n,\varepsilon})\leq Vol(K)^{n-1}%
p_{V^{2},|(\alpha,\beta)|}^{n}(H_{\varepsilon}).
\]
This last inequality implies that the series $\sum_{n\geq1}%
\frac{L_{n,\varepsilon}}{n!}$ normally converges, for the usual
topology of uniform convergence of all the derivatives on any
compact subset of $X^{2}$. Set
$$S_{\varepsilon}=\sum_{n=1}^{+\infty}%
\frac{L_{n,\varepsilon}}{n!}.$$ As $L_{n}$ is in
$\mathcal{G}_{C}(X^{2})$ and since the convergence is uniform,
$S_{\varepsilon}$ belongs to $\mathcal{C}^{\infty}(X^{2})$, for
all $\varepsilon\in(0,1]$. Furthermore, for all compact subset of
$X^{2}$ of the
form $K_{1}\times K_{2}$ and $\left(  \alpha,\beta\right)  \in\mathbb{N}%
^{d}\times\mathbb{N}^{d}$, one has
\begin{align*}
p_{K_{1}\times K_{2},|(\alpha,\beta)|}(S_{\varepsilon})  &  \leq\sum
_{n=1}^{+\infty}\frac{1}{n!}p_{K_{1}\times K_{2},|(\alpha,\beta)|}%
(L_{n,\varepsilon})\\
&  \leq\sum_{n=1}^{+\infty}\frac{1}{n!}Vol(K)^{n-1}p_{V^{2},|(\alpha,\beta
)|}^{n}(H_{\varepsilon})\\
&  \leq\frac{1}{Vol(K)}\left[  e^{Vol(K)p_{V^{2},|(\alpha,\beta)|}%
(H_{\varepsilon})}-1\right]  .
\end{align*}
Since $H$ is in $\mathcal{G}_{C\,ln}(X^{2})$, $p_{V^{2},|(\alpha,\beta
)|}(H_{\varepsilon})=O(|\ln\varepsilon|)$ as $\varepsilon\rightarrow0$, that
is there exists $k\in\mathbb{N}$ such that $p_{V^{2},|(\alpha,\beta
)|}(H_{\varepsilon})\leq\ln\left(  \frac{1}{\varepsilon^{k}}\right)  $, so
\[
p_{K_{1}\times K_{2},|(\alpha,\beta)|}(S_{\varepsilon})\leq C_{K}%
\varepsilon^{-k Vol(K)},
\]
where $C_{K}$ is a constant depending only on $K$ and not on the
representative of $H$. Consequently, $(S_{\varepsilon})$ is in $\mathcal{E}%
_{M}(X^{2})$ and we denote by $S$ its class in $\mathcal{G}(X^{2})$.

The independence of $S$ with respect to the representatives is
classically proved by taking two representatives of $H$, which
gives two sums, $\left( S_{\varepsilon}^{1}\right)  $ and $\left(
S_{\varepsilon}^{2}\right)$, obtained by the process described
above, and by estimating the difference $\left(
S_{\varepsilon}^{1}-S_{\varepsilon}^{2}\right)$. This uses similar
estimates as above.\ Finally, the assertion concerning the support
is proved with similar arguments as the ones of the proof of
theorem \ref{corol}.\medskip

In \cite{BCD1}, it is shown that the exponential defined by
theorem \ref{defexp1} inherits the main expected functional
properties.

\begin{proposition}
If $H$ is in $\mathcal{G}_{Cln}(X^{2})$ then%
\[
\widehat{H}\circ e^{\widehat{H}}=e^{\widehat{H}}\circ\widehat{H}%
\ ;\ \ \ \ e^{a\widehat{H}}\circ e^{b\widehat{H}}=e^{(a+b)\widehat{H}%
},\ \text{for all }\left(  a,b\right)  \in\mathbb{R}^{2}\ ;\ \ \ \frac{d}%
{dt}~~e^{t\widehat{H}}=\widehat{H}\circ e^{t\widehat{H}}.
\]

\end{proposition}
By applying theorem \ref{caracker} concerning the characterization
of generalized integral operators by their kernel, these
properties are proved by using the associated kernels.

\end{document}